\documentclass[10pt]{article}
\usepackage{amsmath} %provare a modificare per toglierli
\usepackage{amssymb}
\usepackage{latexsym}
\usepackage{xcolor}
\usepackage{fancyhdr}
\allowdisplaybreaks 
\usepackage{todonotes}

 \usepackage{comment}

\def\XXint#1#2#3{{\setbox0=\hbox{$#1{#2#3}{\int}$} 
\vcenter{\hbox{$#2#3$}}\kern-.5\wd0}}   
 \numberwithin{equation}{section}
\newtheorem{theorem}[equation]{Theorem}
\newtheorem{proposition}[equation]{Proposition}

\newtheorem{remark}[equation]{Remark}
\newtheorem{lemma}[equation]{Lemma}

\title{Regularity properties of certain  convolution operators in H\"{o}lder spaces}

 \author{  
Matteo Dalla Riva\thanks{ Dipartimento di Tecnica e Gestione dei Sistemi Industriali, 
Universit\`a degli Studi di Padova, 
Stradella San Nicola  3, Vicenza 36100, Italy. 
  Email: matteo.dallariva@unipd.it}, Massimo Lanza de Cristoforis\thanks{
  Dipartimento di Matematica `Tullio Levi-Civita', 
Universit\`a degli Studi di Padova, 
Via Trieste 63, Padova 35121, 
Italy. E-mail: mldc@math.unipd.it} 
and Paolo Musolino\thanks{Dipartimento di Matematica `Tullio Levi-Civita', 
Universit\`a degli Studi di Padova, 
Via Trieste 63, Padova 35121, 
Italy. E-mail: paolo.musolino@unipd.it}
}

\date{\ }
 
\begin{document}

 \maketitle

\noindent
{\bf Abstract:}  The aim of this paper is to prove a theorem of C.~Miranda on the H\"older regularity of convolution operators acting on the boundary of an open set in the limiting case in which the open set is of class $C^{1,1}$ and  the densities are of class $C^{0,1}$. The convolution operators that we consider are generalizations of those that are associated to layer potential operators, which are a useful tool for the analysis of boundary value problems. 

 \vspace{\baselineskip}



\section{Introduction}  Potential theory is a powerful method to study boundary value problems arising in several real-world problems. Here, as an example we mention applications to inverse problems (cf.~Kirsch \cite{Ki21}), to electromagnetism (Angell and Kirsch \cite{AnKi04}, Kirsch and Hettlich \cite{KiHe15}), to fluid mechanics (Kohr and Pop \cite{KoPo04}), and to imaging (Ammari {\it et al.}~\cite{AmBrGaKaLeWa15}). By potential theory, one can convert a boundary value problem for an elliptic or a parabolic partial differential equation into a system of integral equations.  As a consequence, in view of the importance of such method, several authors have devoted monographs to the study of integral operators arising in potential theory in various settings. We mention the  book by McLean \cite{Mc00} on potential theory in Lipschitz domains, the monograph by Hsiao and Wendland \cite{HsWe21}, the series by Mitrea, Mitrea, and Mitrea \cite{MiMiMi22a, MiMiMi22b, MiMiMi23a, MiMiMi23b, MiMiMi23c},  the work by Medkov\'a \cite{Me18}, and  \cite{DaLaMu21}. Of course, if one is interested into investigating boundary value problems by potential theory in a specific framework, then one needs to study mapping properties of the involved integral operators ({\it e.g.}, the single and the double layer operators) in the corresponding framework, such as the one of Sobolev spaces or that of   H\"{o}lder and  Schauder spaces. Such operators appear as convolution operators acting on the boundary of an open set, and specific mapping properties of such type of functions is the topic of the present article. 

In a remarkable paper {providing a general treatment for the H\"older regularity of potential operators, Miranda \cite{Mi65}  has shown that 
if $\alpha\in]0,1[$, $m$, $n\in {\mathbb{N}}$, $m\geq 1$, $n\geq 2$, $\Omega$ is a bounded open subset of ${\mathbb{R}}^n$ of class $C^{m,\alpha}$, $k$ is an odd positively homogeneous function of degree $-(n-1)$ which is of class $C^{2m}$, 		$\mu$ is  a	function on the boundary $\partial\Omega$ of $\Omega$  of class $C^{m-1,\alpha}$,   and
\[
K[k,\mu](x)\equiv\int_{\partial\Omega}k(x-y)\mu(y)\,d\sigma_y
\qquad\forall x\in {\mathbb{R}}^n\setminus\partial\Omega\,,
\] 
then $K[k,\mu]_{|\Omega}$ extends to  a	  function   $K[k,\mu]^+$ on the closure $\overline{\Omega}$ of $\Omega$
  of class $C^{m-1,\alpha}$.  
Moreover, Miranda has proved that the map from $ C^{m-1,\alpha}(\partial\Omega)$ to $C^{m-1,\alpha}(\overline{\Omega})$ that  takes $\mu$ to $K[k,\mu]^+$ is linear and continuous and that a corresponding statement holds in the exterior of $\Omega$. The result of Miranda \cite{Mi65} generalizes previous results in the half space of Agmon Douglis and Nirenberg \cite{AgDoNi59} and under stronger assumptions on the boundary of $\Omega$  of Avantaggiati \cite{Av63}.  For further developments of the ideas contained in Miranda's proof, see, for example, Cialdea \cite{Ci95} and Cialdea, Leonessa, and Malaspina \cite{CiLeMa19}.	 
We also mention  that Wiegner~\cite{Wi93} has proved that if  a multi-index  $\gamma\in {\mathbb{N}}^{n}$ has odd length and $\Omega$ is of class $C^{m,\alpha}$, then the operator with kernel 
\[
 (x-y)^{\gamma}  |x-y|^{-(n-1)- |\gamma|}  
\]
is continuous from  $C^{m-1,\alpha}(\partial\Omega)$ to $C^{m-1,\alpha}(\overline{\Omega})$ (and a corresponding result for the exterior of $\Omega$) and the characterization of Mitrea, Mitrea and Verdera \cite[Theorem 1]{MitMitVe16} of the regularity of $\Omega$ by means of the continuity of the Riesz
transforms on $\partial \Omega$.

In \cite[Theorem 4.17]{DaLaMu21}, the authors have simplified the proof of Miranda \cite{Mi65} and weakened the assumptions on $k$ and proved  a continuity theorem in both the variables $k$ and $\mu$ in case $m=1$. In this paper, we plan to consider the theorem of Miranda in case $m=1$ and $\alpha=1$ and show that if  $\Omega$ is of class $C^{1,1}$, then 
 the map from the cartesian product of the space of odd positively homogeneous functions of degree $-(n-1)$ and locally of class $C^{1,1}$ times $C^{0,1}(\partial\Omega)$ to $C^{0,\omega_1 }(\overline{\Omega})$
  which takes the pair $(k,\mu)$ to $K^+[k,\mu]$ is bilinear and continuous, and that a corresponding statement holds in the exterior of $\Omega$ (see Theorem \ref{thm:mirandao01}).
  Here $C^{0,\omega_{1} }(\overline{\Omega})$ denotes the space
  of functions  that satisfy a generalized $\omega_{1}$-H\"{o}lder condition with
\[
\omega_{1}(r)
\equiv
\left\{
\begin{array}{ll}
0 &r=0\,,
\\
r |\ln r | &r\in]0, e^{-1}]\,,
\\
 e^{-1}  & r\in ] e^{-1},+\infty[\,,
\end{array}
\right.
\]
(see equation (\ref{omth})  below	 	with
$
r_{1}\equiv e^{-1}
$). For the classical definition of the generalized H\"{o}lder or Schauder spaces we refer the reader to
  \cite[p.~75, \S 2]{DoLa17} and to \cite[\S 2.6, \S 2.11]{DaLaMu21}.

\section{Preliminaries and notation}\label{sec:tecprel}
 
If $A$  is a matrix with real (or complex) entries, 
       $A^{t}$ denotes the transpose matrix of $A$. The symbol $O_n({\mathbb{R}})$ denotes the orthogonal group.  
       The symbol
$| \cdot|$ denotes the Euclidean modulus   in
${\mathbb{R}}^{n}$ or in ${\mathbb{C}}$. For all $r\in]0,+\infty[$, $ x\in{\mathbb{R}}^{n}$, 
$x_{j}$ denotes the $j$-th coordinate of $x$, and  
 ${\mathbb{B}}_{n}( x,r)$ denotes the ball $\{y\in{\mathbb{R}}^{n}:\, |x- y|<r\}$.  If ${\mathbb{D}}$ is a subset of $ {\mathbb{R}}^n$, 
then we set
\[
B({\mathbb{D}})\equiv\left\{
f\in {\mathbb{C}}^{\mathbb{D}}:\,f\ \text{is\ bounded}
\right\}
\,,\quad
\|f\|_{B({\mathbb{D}})}\equiv\sup_{\mathbb{D}}|f|\qquad\forall f\in B({\mathbb{D}})\,.
\]
Then $C^0({\mathbb{D}})$ denotes the set of continuous functions from ${\mathbb{D}}$ to ${\mathbb{C}}$ and we introduce the subspace
$
C^0_b({\mathbb{D}})\equiv C^0({\mathbb{D}})\cap B({\mathbb{D}})
$
of $B({\mathbb{D}})$.  Let $\omega$ be a function from $[0,+\infty[$ to itself such that
\begin{eqnarray}
\nonumber
&&\qquad\qquad\omega(0)=0,\qquad \omega(r)>0\qquad\forall r\in]0,+\infty[\,,
\\
\label{om}
&&\qquad\qquad\omega\ {\text{is\   increasing,}}\ \lim_{r\to 0^{+}}\omega(r)=0\,,
\\
\nonumber
&&\qquad\qquad{\text{and}}\ \sup_{(a,t)\in[1,+\infty[\times]0,+\infty[}
\frac{\omega(at)}{a\omega(t)}<+\infty\,.
\end{eqnarray}
Here `$\omega$ is increasing' means that 
$\omega(r_1)\leq \omega(r_2)$ whenever $r_1$, $r_2\in [0,+\infty[$ and $r_1<r_2$.
If $f$ is a function from a subset ${\mathbb{D}}$ of ${\mathbb{R}}^n$   to ${\mathbb{C}}$,  then we denote by   $|f:{\mathbb{D}}|_{\omega  }$  the $\omega$-H\"older constant  of $f$, which is delivered by the formula   
\[
%\label{om1a}
|f:{\mathbb{D}}|_{\omega  
}
\equiv
\sup\left\{
\frac{|f( x )-f( y)|}{\omega(| x- y|)
}: x, y\in {\mathbb{D}} ,  x\neq
 y\right\}\,.
\]        
If $|f:{\mathbb{D}}|_{\omega }<+		\infty$, we say that $f$ is $\omega $-H\"{o}lder continuous. Sometimes, we simply write $|f|_{\omega }$  
instead of $|f:{\mathbb{D}}|_{\omega }$. The
subset of $C^{0}({\mathbb{D}} ) $  whose
functions  are
$\omega $-H\"{o}lder continuous    is denoted  by  $C^{0,\omega } ({\mathbb{D}})$
and $|f:{\mathbb{D}}|_{\omega }$ is a semi-norm on $C^{0,\omega } ({\mathbb{D}})$.  
Then we consider the space  $C^{0,\omega }_{b}({\mathbb{D}} ) \equiv C^{0,\omega } ({\mathbb{D}} )\cap B({\mathbb{D}} ) $ with the norm \[
%\label{om1b}
\|f\|_{ C^{0,\omega }_{b}({\mathbb{D}} ) }\equiv \sup_{x\in {\mathbb{D}} }|f(x)|+|f|_{\omega }\qquad\forall f\in C^{0,\omega }_{b}({\mathbb{D}} )\,.
\] 
\begin{remark}
\label{rem:om4}
Let $\omega$ be as in (\ref{om}). 
Let ${\mathbb{D}}$ be a   subset of ${\mathbb{R}}^{n}$. Let $f$ be a bounded function from $ {\mathbb{D}}$ to ${\mathbb{C}}$, $a\in]0,+\infty[$.  Then,
\[
\label{rem:om5}
\sup_{x,y\in {\mathbb{D}},\ |x-y|\geq a}\frac{|f(x)-f(y)|}{\omega(|x-y|)}
\leq \frac{2}{\omega(a)} \sup_{{\mathbb{D}}}|f|\,.
\]
\end{remark}
 In the case in which $\omega $ is the function 
$r^{\alpha}$ for some fixed $\alpha\in]0,1]$, a so-called H\"{o}lder exponent, we simply write $|\cdot:{\mathbb{D}}|_{\alpha}$ instead of
$|\cdot:{\mathbb{D}}|_{r^{\alpha}}$, $C^{0,\alpha} ({\mathbb{D}})$ instead of $C^{0,r^{\alpha}} ({\mathbb{D}})$, $C^{0,\alpha}_{b}({\mathbb{D}})$ instead of $C^{0,r^{\alpha}}_{b} ({\mathbb{D}})$, and we say that $f$ is $\alpha$-H\"{o}lder continuous provided that 
$|f:{\mathbb{D}}|_{\alpha}<	+	\infty$. For each $\theta\in]0,1]$, we define the function $\omega_{\theta}(\cdot)$ from $[0,+\infty[$ to itself by setting
\begin{equation}
\label{omth}
\omega_{\theta}(r)\equiv
\left\{
\begin{array}{ll}
0 &r=0\,,
\\
r^{\theta}|\ln r | &r\in]0,r_{\theta}]\,,
\\
r_{\theta}^{\theta}|\ln r_{\theta} | & r\in ]r_{\theta},+\infty[\,,
\end{array}
\right.
\end{equation}
where
$
%\label{omth1}
r_{\theta}\equiv e^{-1/\theta}
$ for all $\theta\in ]0,1]$. Obviously, $\omega_{\theta}  $ is concave and satisfies   condition (\ref{om}).
 We also note that if ${\mathbb{D}}\subseteq {\mathbb{R}}^n$, then the continuous embedding
\begin{equation}\label{eq:0t0om0t}
C^{0, \theta }_b({\mathbb{D}})\subseteq 
C^{0,\omega_\theta }_b({\mathbb{D}})\subseteq 
C^{0,\theta'}_b({\mathbb{D}})
\end{equation}
holds for all $\theta'\in ]0,\theta[$. For the standard properties of the spaces of H\"{o}lder or Lipschitz continuous functions, we refer    to \cite[\S 2]{DoLa17}, \cite[\S 2.6]{DaLaMu21}.\par

Let $\Omega$ be an open subset of ${\mathbb{R}}^n$. The space of $m$ times continuously 
differentiable real-valued functions on $\Omega$ is denoted by 
$C^{m}(\Omega,{\mathbb{C}})$, or more simply by $C^{m}(\Omega)$. 
Let $f\in C^{m}(\Omega) $. Then   $Df$ denotes the Jacobian matrix of $f$. 
Let  $\eta\equiv
(\eta_{1},\dots ,\eta_{n})\in{\mathbb{N}}^{n}$, $|\eta |\equiv
\eta_{1}+\dots +\eta_{n}  $. Then $D^{\eta} f$ denotes
$\frac{\partial^{|\eta|}f}{\partial
x_{1}^{\eta_{1}}\dots\partial x_{n}^{\eta_{n}}}$.    The
subspace of $C^{m}(\Omega )$ of those functions $f$ whose derivatives $D^{\eta }f$ of
order $|\eta |\leq m$ can be extended with continuity to 
$\overline{\Omega}$  is  denoted $C^{m}(
\overline{\Omega})$. 
The
subspace of $C^{m}(\overline{\Omega} ) $  whose
functions have $m$-th order derivatives that are
H\"{o}lder continuous  with exponent  $\alpha\in
]0,1]$ is denoted $C^{m,\alpha} (\overline{\Omega})$ and the
subspace of $C^{m}(\overline{\Omega} ) $  whose
functions have $m$-th order derivatives that are
$\omega $-H\"{o}lder continuous    is denoted $C^{m,\omega } (\overline{\Omega})$. 

 Then $C^{m,\omega }_{{\mathrm{loc}}}(\overline{\Omega }) $  denotes 
the space  of those functions $f\in C^{m}(\overline{\Omega} ) $ such that $f_{|\overline{\Omega'}} $ belongs to $
C^{m,\omega }(   \overline{ \Omega' }   )$ for all bounded open subsets $\Omega'$ of ${\mathbb{R}}^n$ such that $\overline{\Omega'}\subseteq\overline{\Omega}$.  
 
Now let $\Omega $ be a bounded
open subset of  ${\mathbb{R}}^{n}$. Then $C^{m}(\overline{\Omega} )$,  $C^{m,\omega  }(\overline{\Omega})$ with $\omega$ as in (\ref{om}),
  $C^{m,\alpha }(\overline{\Omega})$,  are equipped   with their usual norm and are well known to be 
Banach spaces  (cf.~\textit{e.g.},  \cite[\S 2]{DoLa17},  \cite[\S 2.11]{DaLaMu21}).  

 For the definition of a bounded open Lipschitz subset of ${\mathbb{R}}^{n}$ and  for  the (classical) definition of open set of class $C^{m}$ or of class $C^{m,\alpha}$, we refer for example to \cite[\S 2.9, \S 2.13]{DaLaMu21}. 
 
  For the (classical) definition of the H\"{o}lder and Schauder spaces     $C^{m,\omega }(\partial\Omega)$, $C^{m,\alpha}(\partial\Omega)$ 
on the boundary $\partial\Omega$ of an open set $\Omega$ for some $m\in{\mathbb{N}}$, $\alpha\in ]0,1]$ $\omega $ as in  (\ref{om}),  we refer for example to  \cite[\S  2.20]{DaLaMu21}, \cite[\S 2]{DoLa17}.	
  
Finally, the symbol 
$m_n$ denotes the
$n$ dimensional Lebesgue measure.

 We now introduce some (known) preliminary statements. 
\begin{lemma}
\label{lem:comin} Let $n\in {\mathbb{N}}\setminus\{0,1\}$.
Let $\Omega$ be a bounded open Lipschitz subset of ${\mathbb{R}}^{n}$. Then the following statements hold.
\begin{enumerate}
\item[(i)]  If $\lambda\in ]-\infty, n-1[$, then\index{$c'_{\Omega,\lambda}$} 
\[
c'_{\Omega,\lambda}\equiv \sup_{x\in\partial\Omega}
\int_{\partial\Omega}\frac{d\sigma_{y}}{|x-y|^{\lambda}}
\]
is finite.
\item[(ii)] If $\lambda\in ]-\infty, n-1[$, then\index{$c''_{\Omega,\lambda}$}
\[
c''_{\Omega,\lambda}
\equiv \sup_{x'\in\partial\Omega,s\in]0,+\infty[}
s^{\lambda    -(n-1)}\int_{{\mathbb{B}}_{n}(x',s)\cap\partial\Omega}
\frac{d\sigma_{y}}{|x'-y|^{\lambda}}
\]
is finite.
\item[(iii)] If  $\lambda\in ]n-1,+\infty[$, then\index{$c'''_{\Omega,\lambda}$}
\[
c'''_{\Omega,\lambda}
\equiv
 \sup_{x'\in\partial\Omega,s\in]0,+\infty[}
 s^{\lambda -(n-1) }
 \int_{\partial\Omega\setminus
 {\mathbb{B}}_{n}(x',s)
 }
 \frac{d\sigma_{y}}{|x'-y|^{\lambda}} 
\]
is finite.
\item[(iv)] 
\[
c^{iv}_{\Omega}
\equiv
 \sup_{x'\in\partial\Omega,\ 
 s\in]0,1/e[
 }
 |\ln s |^{-1}
 \int_{\partial\Omega\setminus
 {\mathbb{B}}_{n}(x',s)
 }
 \frac{d\sigma_{y}}{|x'-y|^{n-1}}<+\infty\,.
\]
\end{enumerate}
\end{lemma}
{\bf Proof.} For a proof of statements (i)--(iii), we refer to \cite[Lemma 2.53]{DaLaMu21}. We now turn to   statement (iv). By \cite[Proposition 3.1]{La24a}, 
\begin{eqnarray*}\nonumber
&&\text{there\ exists}\  a\in]0,1/(2e)[\ \text{such\ that\  for\ each}\  \rho\in]0,a[\ \text{and}\  x'\in \partial\Omega,
 \\ \label{prop:mfbd0}
&&\text{there\ exists}\  x''\in \partial\Omega\ \text{such\ that}\  |x'-x''|=\rho\,.
\end{eqnarray*}
Then \cite[Lemma 3.5 (iv)]{DoLa17} implies that the supremum in (iv) is finite if we take $s\in ]0,1/(2e)[$. Since
\[
|\ln s |^{-1}
 \int_{\partial\Omega\setminus
 {\mathbb{B}}_{n}(x',s)
 }
 \frac{d\sigma_{y}}{|x'-y|^{n-1}}
 \leq
 |\ln (1/(2e)) |^{-1}  \frac{m_{n-1}(\partial\Omega)}{(1/(2e))^{n-1}}
 \quad\forall s\in[1/(2e),1/e[
 \,,
 \]
statement (iv) holds true.\hfill  $\Box$ 

\vspace{\baselineskip}

\section[A Sufficient Condition for  $\omega_1 $-H\"{o}lder Continuity]{A sufficient condition for the $\omega_1 $-H\"{o}lder continuity of continuously differentiable functions}\label{sec:holdercond}

This section aims at presenting a sufficient condition for the $\omega_1 $-H\"{o}lder continuity of continuously differentiable functions, building on the work of Carlo Miranda.

We first introduce the following elementary lemma.
\begin{lemma}\label{lem:omth}
 The function $\omega_1$ of (\ref{omth}) satisfies the following inequality
 \[
| \omega_1(t'')-\omega_1(t')|\leq \omega_1(t''-t')\qquad\forall t',t''\in [0,+\infty
 [\,.	  
 \]
 In particular, $\omega_1$ belongs to $C^{0,\omega_1 }([0,+\infty[)$. 
\end{lemma}
{\bf Proof.} Let $t',t''\in [0,+\infty[$, $t'\neq t''$. There is no loss of generality in assuming that $t'<t''$. Since $\omega_1$ is strictly increasing in $[0,1/e]$ and constant on $[1/e,+\infty[$, we  
  estimate $\omega_1(t'')-\omega_1(t')$ by considering some separate cases depending on whether $t'$ and $t''$ belong to either $[0,1/e]$ or $[1/e,+\infty[$.  
  
  We first consider the case in which $t',t''\in ]0,1/e]$. Since $\omega_1$ is strictly increasing in $]0,1/e]$ and $-\log$ is strictly decreasing in $]0,1/e]$ and $t''-t'\leq t''$, we have
  \begin{eqnarray*}
\lefteqn{
0<\omega_1(t'')-\omega_1(t')= t''|\log t''|-  t'|\log t'|
}
\\ \nonumber
&&\qquad
=(t''-t')(-\log t'')+t'(-\log t''+\log t')
\\ \nonumber
&&\qquad
=(t''-t')(-\log t'')-t'(\log (t''/t'))
\\ \nonumber
&&\qquad
\leq(t''-t')(-\log (t''-t'))+0=\omega_{1}(t''-t')\,.
\end{eqnarray*}
Hence, we deduce the validity of the inequality 
 \begin{equation}\label{lem:omth1}
 0<\omega_1(t'')-\omega_1(t')\leq \omega_{1}(t''-t')
 \end{equation}
  for all $t'$, $t''\in [0,1/e]$, $t'<t''$. Next we consider   the case in which
  $t'\leq 1/e\leq t''$. Since $\omega_1(t'')=\omega_1(1/e)$ and $\omega_1$ is increasing, inequality (\ref{lem:omth1})   implies that 
    \[
  0\leq\omega_1(t'')-\omega_1(t')=\omega_1(1/e)-\omega_1(t')\leq 
  \omega_{1}(1/e-t')\leq \omega_{1}(t''-t')\,.
  \]
 Finally,  we consider   the case in which
  $1/e\leq t'<t''$. Since $\omega_1$ is constant in $[1/e,+\infty[$, we have
  \[
  0\leq\omega_1(t'')-\omega_1(t')=0\leq \omega_{1}(t''-t')
  \]
and thus the proof is complete.\hfill  $\Box$ 

\vspace{\baselineskip}

To carry out our proofs, we need a   sufficient condition for a real valued continuously differentiable function $f$ in an open subset $\Omega$ of ${\mathbb{R}}^n$  to be $\omega_1$-H\"{o}lder continuous in $\Omega$ that generalizes a corresponding sufficient condition for the $\alpha$-H\"{o}lder continuity of $f$ in $\Omega$  for some  $\alpha\in ]0,1]$ 
of \cite[\S 2.19]{DaLaMu21}. Notice that a $C^1$ function on $\Omega$ is always $\omega_1$-H\"{o}lder continuous on any compact subset of $\Omega$, but it may fail to remain so  in a neighborhood of  the boundary $\partial\Omega$.

Our generalization applies to sets of class $C^1$. To do so, we need to introduce an appropriate neighborhood of $\partial\Omega$ in ${\mathbb{R}}^n$ and we want such a  neighborhood to be `globally parametrized.' Thus it is natural to think of the 
 set
 \[
\{x+t\nu_\Omega(x):\,t\in]-t_0,t_0[\,,\ x\in\partial\Omega\}\,,
\]
for some small $t_0\in]0,+\infty[$. 
Unfortunately however, the map $\nu_\Omega$ is only of class 
$C^{0}$  if $\Omega$ of class $C^{1}$.  	Thus, the corresponding `global correspondence' 
\[
(t,x)\mapsto x+t\nu_\Omega(x)
\]
between the points of $]-t_0,t_0[\times (\partial\Omega)$ and the points of the above set cannot be expected to be Lipschitz continuous  as we need it to be  to prove that such a correspondence is injective. Then  we follow the paper \cite{LaRo08} with Rossi  and we state  the following two lemmas  (cf.~\cite[Lemmas 2.78, 2.79]{DaLaMu21}), which show  the existence of a sufficiently small $t_0\in]0,+\infty[$ and of
a unit vector field $a$ of class $C^{\infty}$ in a neighborhood of $\partial\Omega$     such that the map from 
$]-t_0,t_0[\times (\partial\Omega)$ to ${\mathbb{R}}^n$ that takes $(t,x)$ to $x+ta(x)$ is injective and Lipschitz continuous   and such that
the set
\[
\{x+ta(x):\,t\in]-t_0,t_0[\,,\ x\in\partial\Omega\}
\]
is an open neighborhood of $\partial\Omega$ in ${\mathbb{R}}^n$.  \begin{lemma}\label{lem:exacu}
 Let $n\in {\mathbb{N}}\setminus\{0,1\}$. Let $\Omega$  be a bounded open subset of ${\mathbb{R}}^n$ of class $C^1$. Let $\vartheta\in]0,1[$. Then there exist an open neighborhood $U$ of class $C^\infty$ of $\partial\Omega$, $a\in C^\infty(\overline{U},{\mathbb{R}}^n)$, and $\tau\in]0,+\infty[$ such that the following conditions are fulfilled
\begin{eqnarray}\label{lem:exacu1}
 &&|a(x)|=1\qquad\forall x\in\overline{U}\,,
 \\ \nonumber
 &&\sup_{\partial\Omega}|a-\nu_\Omega|<\vartheta\,,
 \\ \nonumber
 &&\inf_{\partial\Omega} a\cdot \nu_\Omega>1-\frac{\vartheta^2}{2}>1-\vartheta\,,
 \\ \nonumber
 && |a(x)\cdot (y-x)|\leq \vartheta |x-y|
 \qquad\forall x,y\in \partial\Omega\,,\ |x-y|< \tau\,.
 \end{eqnarray}
 \end{lemma}
We note that Fichera~\cite[pp.~207-208]{Fi55}  already recognized the importance of the existence of a vector field $a$ such that the essential infimum of $a \cdot\nu_\Omega  $ is bounded away from $0$  for sets with a piecewise smooth boundary in the analysis of boundary value problems for systems of partial differential equations, and in particular for the analysis of boundary value problems with unilateral constraints (see also Fichera~\cite[p.~413]{Fi84}). For the introduction of a vector field such that the essential infimum of $a \cdot\nu_\Omega  $ is bounded away from $0$ in the case of Lipschitz sets, we mention  Grisvard~\cite[Lemma~1.5.1.9]{Gr85}.\par

 A unit vector field $a$ as in the previous lemma can play the role of a normal in the definition of a tubular neighborhood of $\partial\Omega$, as the following statement shows  (see also reference \cite{LaRo08} with Rossi for a related result which holds under stronger assumptions on $\Omega$).
\begin{lemma}\label{lem:tubnb}
Let $n\in {\mathbb{N}}\setminus\{0,1\}$. Let $m\in{\mathbb{N}}\setminus\{0\}$. Let $\Omega$ be a bounded open subset of ${\mathbb{R}}^n$ of class $C^1$.   If 
an open neighborhood $U$ of  $\partial\Omega$ 
and $a\in C^1(\overline{U},{\mathbb{R}}^n)$  satisfy conditions  (\ref{lem:exacu1}) for some $\vartheta\in]0,1[$ and $\tau\in ]0,+\infty[$, then there exists $t_1\in]0,+\infty[$ such that the following statements hold.
\begin{enumerate}
\item[(i)] The map $\Psi$ from $(\partial\Omega)\times ]-t_1,t_1[$ to ${\mathbb{R}}^n$ defined by 
\[
\Psi (x,t)=x+t a(x)\qquad \forall (x,t)\in (\partial\Omega)\times ]-t_1,t_1[
\]
is injective, and the set
\[
A(t_2)\equiv \Psi ((\partial\Omega)\times ]-t_2,t_2[)
\]
is an open neighborhood of $\partial\Omega$ for all $t_2\in ]0,t_1]$.
\item[(ii)] We have
\[
x+t a(x)\in \Omega\qquad\forall t\in ]-t_1,0[ 
\quad and\quad
x+t a(x)\in  \mathbb{R}^n\setminus \overline{\Omega} \qquad\forall t\in ]0,t_1[\,.
\]
\end{enumerate}
\end{lemma}
Next we point out the validity of the following elementary lemma (cf. \cite[Lem.~2.80]{DaLaMu21}). 
\begin{lemma}\label{lem:bebelbd} 
 Let $n\in {\mathbb{N}}\setminus\{0,1\}$. Let $\vartheta\in]0,1[$. If $v$, $w\in {\mathbb{R}}^n $ and if
 \[
 \left| v \cdot w \right| \leq \vartheta |v|\,|w|\,,
 \]
 then $|v+w|^2\geq (1-\vartheta)(|v|^2+|w|^2)+\vartheta (|v|-|w|)^2$.
 \end{lemma}

 We say that an open subset $\Omega$ of ${\mathbb{R}}^n$  is of class $C^0$ provided that for each point $p\in\partial\Omega$ there exist
$R_p\in O_{n}({\mathbb{R}})$, $r$, $\delta\in ]0,+\infty[$ such that the intersection
 \[
 R_p(\Omega-p) \cap ({\mathbb{B}}_{n-1}(0,r)\times]-\delta,\delta[    )
 \]
 is the strict hypograph of a continuous function $\gamma_p$ from ${\mathbb{B}}_{n-1}(0,r)$ to $]-\delta ,\delta [$ which vanishes at $0$ and such that $|\gamma_p(\eta)|<\delta/2$ for all $\eta\in {\mathbb{B}}_{n-1}(0,r)$,  i.e., provided that there exists $\gamma_p\in C^{0}({\mathbb{B}}_{n-1}(0,r), ]-\delta ,\delta [)$ such that 
 \begin{eqnarray*}
\label{prelim.cocylind1}
\lefteqn{R_p(\Omega-p )\cap ({\mathbb{B}}_{n-1}(0,r)\times ]-\delta,\delta[) 
}
\\ \nonumber
&&\qquad 
=\left\{
(\eta,y)\in 
{\mathbb{B}}_{n-1}(0,r)\times ]-\delta,\delta[:\, y<\gamma_p(\eta)
\right\}
\equiv{\mathrm{hypograph}}_{s}(\gamma_p) 
\,, 
\\ \nonumber
&&\qquad 
|\gamma_p(\eta)|<\delta/2\qquad\forall \eta\in {\mathbb{B}}_{n-1}(0,r)\,,\qquad
\gamma_p(0)=0\,.
\end{eqnarray*}
Here the subscript `$s$' of  `$\,{\mathrm{hypograph}}_{s}$' stands for `strict'. Then we say that the set
 \[
 C(p,R_p,r,\delta)\equiv p+ R_p^{t}({\mathbb{B}}_{n-1}(0,r)\times]-\delta,\delta[    )
 \]
 is a coordinate cylinder  for $\Omega$ around $p$, that the function $\gamma_p$ represents $\partial\Omega$ in $C(p,R_p,r,\delta)$ as a graph and that the function $\psi_{p}$ from ${\mathbb{B}}_{n-1}(0,r)$ to ${\mathbb{R}}^{n}$ defined by 
\[ 
\psi_{p}(\eta)\equiv p+R_p^{t}\left(
\begin{array}{c}
\eta
\\
\gamma_p(\eta)
\end{array}
\right)\qquad\forall \eta\in {\mathbb{B}}_{n-1}(0,r)
\] 
 is the parametrization of $\partial\Omega$ around $p$ in the coordinate cylinder $C(p,R_p,r,\delta)$.  Then we say that $\Omega$ is of class $C^{m}$ or of class $C^{m,\alpha}$
for some $m\in{\mathbb{N}}$, $\alpha\in ]0,1]$ provided that  $\gamma_p$ is of class $C^{m}$ or of class $C^{m,\alpha}$ for all $p\in\partial\Omega$ (as in \cite[\S 2.7, \S 2.13]{DaLaMu21}). 	Then we show the following variant of \cite[Lem.~2.81]{DaLaMu21}.
\begin{lemma}\label{lem:cylcoomest}
 Let $n\in {\mathbb{N}}\setminus\{0,1\}$. Let $\Omega$  be a bounded open subset of ${\mathbb{R}}^n$ of class $C^1$. Let $\vartheta\in]0,1[$. Let $U$ be an open neighborhood of $\partial\Omega$  of class $C^\infty$.
 Let   $a$ in $ C^\infty(\overline{U},{\mathbb{R}}^n)$
 satisfy the conditions in (\ref{lem:exacu1}) for some $\vartheta\in]0,1[$ and $\tau\in]0,+\infty[$.
  Let $t_1\in]0,1/e[$ be as in Lemma \ref{lem:tubnb}. Let $p\in\partial\Omega$. Let $r$, $\delta\in ]0,+\infty[$, 
\[
2\delta<\min\{\tau/2,t_1\}\,,\qquad r<\delta\,,
\]
 and $R_p\in O_n({\mathbb{R}})$ be such that  $C(p,R_p,r,\delta)$ is a coordinate cylinder  for $\Omega$ around $p$.   Let $\gamma_p\in C^1(\overline{{\mathbb{B}}_{n-1}(0,r)})$ represent $\partial\Omega$ in $C(p,R_p,r,\delta)$ as a graph. Let 
 \[
C(p,R_p,r,\delta)\cap\Omega\subseteq A^+(t_1)\equiv\{ x+t a(x):\,t\in  ]-t_1,0[,x\in\partial\Omega \}\,.
 \]
 Let
 \[
 t_2\in \biggr]0,\min\biggl\{ \frac{r}{4},   \frac{ (1-\vartheta)^{1/2}}{2\sqrt{2}({\mathrm{Lip}}(a)+1)}, \frac{t_1}{2}\biggr\} \biggl[\,.
 \]
 Then there exists $B\in]0,+\infty[$ such that
\begin{equation}\label{lem:cylcoomest1}
|f:\, [p+R_p^t({\mathbb{B}}_{n-1}(0,r/4)\times]-\delta/4,\delta/4[)] \cap A^+(t_2)|_{\omega_1}\leq BM_{t_1, \omega_1	}(f)
\end{equation}
for all $f\in C^1(\Omega)$ such that 
\begin{equation}\label{lem:cylcoomest2}
M_{t_1, \omega_1	 	}(f)\equiv\sup\{
|\log | t|\,|^{-1}|\nabla f(x+ta(x))|:\,(x,t)\in (\partial\Omega)\times ]-t_1,0[ 
\}<+\infty\,,
\end{equation}
(cf.~(\ref{omth}) for the definition of $\omega_1$).
\end{lemma}
{\bf Proof.} Let $p'$, $p''\in [p+R_p^t({\mathbb{B}}_{n-1}(0,r/4)\times]-\delta/4,\delta/4[)]\cap A^+(t_2)$. We plan to estimate the difference $|f(p')-f(p'')|$ and show that it is smaller than a constant times $M_{t_1, \omega_1	}(f)\omega_1(|p'-p''|)$.  In order to exploit condition (\ref{lem:cylcoomest2}) on $f$, we plan 
to 	define a $C^1$ path in $ A^+(t_1)$   that joins $p'$ and $p''$. Lemma \ref{lem:tubnb} implies that there exist   $x'$, $x''\in\partial\Omega$, $t'$, $t''\in ]-t_2,0[$ such that
\[
p'=x'+t'a(x')\,,
\qquad
p''=x''+t''a(x'')\,.
\]
By the triangular inequality, by the equality $|a(x')|=1$, and by the inequality $t_2< r/4 < \delta/4$, we have
\begin{eqnarray*}
\lefteqn{
x'\in {\mathbb{B}}_n(p',t_2)\subseteq [p+R_p^t({\mathbb{B}}_{n-1}(0,r/4)\times]-\delta/4,\delta/4[)]+{\mathbb{B}}_n(0,t_2)
}
\\ \nonumber
&&\qquad\qquad\qquad\qquad\qquad
\subseteq [p+R_p^t({\mathbb{B}}_{n-1}(0,r/2)\times]-\delta/2,\delta/2[)]
\end{eqnarray*} 
and similarly
\[
x''\in [p+R_p^t({\mathbb{B}}_{n-1}(0,r/2)\times]-\delta/2,\delta/2[)]\,.
\]
Next we set
\[
\psi_{p}(\eta)\equiv p+R_p^{t} ( \eta,\gamma_p(\eta))^t
\qquad\forall \eta\in {\mathbb{B}}_{n-1}(0,r)\,.
\] 
Since $x'$, $x''\in\partial\Omega$, then there exist $\xi'$, $\xi''\in {\mathbb{B}}_{n-1}(0,r/2)$ such that
\begin{eqnarray*}
x'&=&\psi_{p}(\xi')=p+R_p^t(\xi',\gamma_p(\xi'))^t\,,
\\
x''&=&\psi_{p}(\xi'')=p+R_p^t(\xi'',\gamma_p(\xi'' ))^t\,.
\end{eqnarray*}
Then we note that the convexity of ${\mathbb{B}}_{n-1}(0,r/2)$ implies that 
\[
y(s)\equiv(1-s)\xi'+s\xi''  \in {\mathbb{B}}_{n-1}(0,r/2)\qquad\forall s\in [0,1]
\]
and that the convexity of $]-t_2,0[$ implies that
\[
(1-s)t'+st''\in ]-t_2,0[\qquad\forall s\in [0,1]\,.
\]
Then we have 
\begin{eqnarray} \nonumber
\lefteqn{G(s)\equiv \psi_p(y(s) )
}
\\ \nonumber
&&\qquad
+[(1-s)t'+st'']a(\psi_p(y(s) ))\in A^+(t_2)\qquad\forall s\in[0,1]\,,
\\ \label{eq:Gsin}
\lefteqn{G(s)\in 
[p+R_p^t({\mathbb{B}}_{n-1}(0,r/2)\times]-\delta/2,\delta/2[)]
+{\mathbb{B}}_{n}(0,t_2)
}
\\ \nonumber
&&\qquad
\subseteq  
[p+R_p^t({\mathbb{B}}_{n-1}(0,3r/4)\times]-3\delta/4,3\delta/4[)]
\quad\forall s\in[0,1]\,, 
\\ \nonumber
\lefteqn{
G(0)=p'\,,\qquad G(1)=p''\,.
}
\end{eqnarray}
Moreover, by arguing as in the proof of  \cite[Lem.~2.81, p.~92]{DaLaMu21}
\[
{\mathrm{length}}(G) \leq |x'-x''|\sqrt{1+{\mathrm{Lip}}^2(\gamma_p)}(1+|t_2|{\mathrm{Lip}}\,(a))+|t'-t''|\,.
\]}
We now wish to estimate $|x'-x''|$ and $|t'-t''|$ in terms of $|p'-p''|$. There is no loss of generality in assuming that
\[
t''<t'\,.
\]
By assumption, we know that
\[
\left|
a(x)\cdot\frac{y-x}{|y-x|}
\right|\leq \vartheta\qquad\forall x,y\in C(p,R_p,r,\delta)\cap(\partial\Omega)\,,\ x\neq y\,.
\]
Indeed, $|x-y|\leq 2r+ \delta<\tau$ for all $ x,y\in C(p,R_p,r,\delta)\cap(\partial\Omega)$. Then Lemma \ref{lem:bebelbd} and the same argument of \cite[Lem.~2.81, (2.72)]{DaLaMu21} imply  that 
\[ 
|t'-t''|\leq \frac{\sqrt{2}}{ (1-\vartheta)^{1/2}}|p'-p''|\,,
\qquad
|x'-x''|\leq \frac{ 2 \sqrt{2} }{ (1-\vartheta)^{1/2}}|p'-p''|
\,.
\] 
Next we take \[
d_1\equiv \min\{
|p'-p''|,r/4
\}
\]
and we define a path with endpoints
\[
x'+(t'-d_1)a(x')\,,
\qquad
x''+(t''-d_1)a(x'')\,.
\]
To do so,  we set
\[
G_{d_1}(s)\equiv G(s)-d_1a(  \psi_p(y(s)) )         \qquad\forall s\in[0,1]\,.
\]
Since
\begin{eqnarray*}
\lefteqn{
(1-s)t'+st''-d_1\in ]-t_2-d_1,0[
}
\\ \nonumber
&&\quad
\subseteq]-(t_1/2)-(r/4),0[\subseteq ]-(t_1/2)-(t_1/8),0[\subseteq]-t_1,0[\qquad\forall s\in]0,1[\,,
\end{eqnarray*}
we have
\[
G_{d_1}(s)\in A^+(t_1) \qquad\forall s\in]0,1[\,.
\]
Then by the membership in \eqref{eq:Gsin} we have 
\begin{eqnarray*}
\lefteqn{G_{d_1}(s)\in  
[p+R_p^t({\mathbb{B}}_{n-1}(0,3r/4)\times]-3\delta/4,3\delta/4[)]+{\mathbb{B}}_n(0,d_1)
}
\\ \nonumber
 &&\qquad\qquad\qquad\qquad\qquad\qquad\qquad\qquad\qquad\qquad
\subseteq C(p,R_p,r,\delta)
\end{eqnarray*}
for all $s\in[0,1]$ and
\[
G_{d_1}(0)=x'+(t'-d_1)a(x')\,,\qquad
G_{d_1}(1)=x''+(t''-d_1)a(x'')\,.
\]
Then, as in \cite[Lem.~2.81, p.~94]{DaLaMu21}, we have the following inequality for the length of $G_{d_1}$,
\begin{eqnarray} \label{prelim.lem:cylcoaest2esa}	 
\lefteqn{
{\mathrm{length}}(G_{d_1})
}
\\ \nonumber
&&\quad  \leq (1+d_1{\mathrm{Lip}}(a))
\\ \nonumber
&&\quad\quad
\times
\left[\frac{2\sqrt{2}}{ (1-\vartheta)^{1/2}}\sqrt{1+{\mathrm{Lip}}^2(\gamma_p) }(1+|t_2|{\mathrm{Lip}}\,(a))+\frac{\sqrt{2}}{ (1-\vartheta)^{1/2}}\right]
 |p'-p''|.
\end{eqnarray}	
 Next we note that
\begin{eqnarray*}
\lefteqn{
x'+(t'-sd_1)a(x')\in [p+R_p^t({\mathbb{B}}_{n-1}(0,3r/4)\times]-3\delta/4,3\delta/4[)]+{\mathbb{B}}_n(0,d_1)
}
\\ \nonumber
&&\qquad\qquad\qquad\qquad\qquad\qquad\qquad\qquad\quad\quad
\subseteq C(p,R_p,r,\delta)
\qquad \forall s\in[0,1]
\end{eqnarray*}
and similarly
\[
x''+(t''-sd_1)a(x'')\in C(p,R_p,r,\delta)
\qquad\forall s\in[0,1]\,.
\]
Then  by the memberships
\begin{eqnarray*}
\lefteqn{
t'-sd_1, t''-sd_1\in ]-t_2-d_1,0[
}
\\ \nonumber
&&\qquad
\subseteq]-(t_1/2)-(r/4),0[\subseteq]-(t_1/2)-(t_1/8),0[\subseteq ]-t_1,0[
\qquad\forall s\in[0,1]\,,
\end{eqnarray*}
we have 
\[
x'+(t'-sd_1)a(x')\in A^+(t_1)\,,\quad      x''+(t''-sd_1)a(x'')\in A^+(t_1)\qquad\forall s\in[0,1]\,.
\]
We now assume that $f\in C^1(\Omega)$  satisfies condition (\ref{lem:cylcoomest2}) and we turn to estimate  $|f(p')-f(p'')|$. We have 
\begin{eqnarray*}
\lefteqn{
|f(p')-f(p'')|=| f(x'+t'a(x'))-  f(x''+t''a(x''))|
}
\\ \nonumber
&&\qquad
\leq | f(x'+t'a(x'))-f(x'+(t'-d_1)a(x'))|
\\ \nonumber
&&\qquad\quad
+| f(x'+(t'-d_1)a(x'))- f(x''+(t''-d_1)a(x'')) |
\\ \nonumber
&&\qquad\quad
+| f(x''+(t''-d_1)a(x''))-f(x''+t''a(x''))|
\\ \nonumber
&&\qquad
\leq\int_{t'-d_1}^{t'}\left|
a(x')\cdot \nabla f(x'+sa(x'))
\right|\,ds
+
\int_0^1|\nabla f(G_{d_1}(s))|\,|G_{d_1}'(s)|\,ds
\\ \nonumber
&&\qquad\quad
+\int_{t''-d_1}^{t''}\left|
a(x'')\cdot\nabla f(x''+sa(x''))
\right|\,ds
\\ \nonumber
&&\qquad
\leq\int_{t'-d_1}^{t'}M_{t_1, \omega_1	}(f)|\log |s|\,|\,ds
\\ \nonumber
&&\qquad\quad
+ \int_0^1M_{t_1, \omega_1	}(f)|\log |(1-s)t'+st''-d_1|\,|\,|G_{d_1}'(s)|\,ds
\\ \nonumber
&&\qquad\quad
+\int_{t''-d_1}^{t''}M_{t_1, \omega_1	}(f)|\log |s|\,|\,ds\,.
\end{eqnarray*}
Since  
\[
(1-s)t'+st''-d_1\in ]-t_2-d_1,-d_1[\qquad\forall s\in[0,1]\,,
\]
we have
\begin{eqnarray}\label{prelim.lem:cylcoaest3}
\lefteqn{
|f(p')-f(p'')|\leq M_{t_1, \omega_1}(f)\int_{|t'|}^{|t'|+d_1}|\log |s|\,|\,ds
}
\\ \nonumber
&&\quad
+{\mathrm{length}}(G_{d_1})|\log |d_1|\,| M_{t_1, \omega_1}(f)
+M_{t_1, \omega_1	}(f) \int_{|t''|}^{|t''|+d_1}|\log |s|\,|\,ds
\\ \nonumber
&&
\leq M_{t_1, \omega_1	}(f)\biggl\{\biggr.
\omega_1(d_1)+d_1
+|p'-p''||\log |d_1|\,|(1+d_1{\mathrm{Lip}}(a))
\\ \nonumber
&&\quad\times
\biggr[\biggl.
\frac{2\sqrt{2} }{ (1-\vartheta)^{1/2}}(1+{\mathrm{Lip}} (\gamma_p))(1+|t_1|{\mathrm{Lip}} (a) )
+\frac{\sqrt{2}}{ (1-\vartheta)^{1/2}}
\biggl.\biggr]
+ \omega_1(d_1)+d_1
\biggl.\biggr\} 
\end{eqnarray}
 (cf. Lemma \ref{lem:omth} and (\ref{prelim.lem:cylcoaest2esa})).  Next,  	we observe that
\[
d_1\leq		\omega_1(d_1)\le\omega_1( |p'-p''|)
\] 
and that 
\begin{eqnarray*}
\lefteqn{
 |p'-p''| |\log |d_1|\,|= |p'-p''| |\log |\min \{|p'-p''|,r/4\}|\,|
}
\\ \nonumber
&& \qquad
\leq\left\{
\begin{array}{ll}
\omega_1( |p'-p''|) &{\text{if}}\ |p'-p''|\leq r/4\,,
 \\
 \frac{{\mathrm{diam}}(\Omega)}{r/4}\omega_1(r/4)
 \leq  \frac{{\mathrm{diam}}(\Omega)}{r/4}\omega_1( |p'-p''|)	&{\text{if}}\ |p'-p''|\geq r/4\,.
\end{array}
\right.
\end{eqnarray*}
Hence, inequality (\ref{prelim.lem:cylcoaest3}) implies the validity of inequality (\ref{lem:cylcoomest1}) and the proof is complete.

\hfill  $\Box$ 

\vspace{\baselineskip}

We are now ready to prove a sufficient condition for the $\omega_1 $-H\"older
  continuity of a continuously differentiable  real valued 
  function defined on an open
 subset of $\mathbb{R}^n$ of class $C^1$
 by modifying an argument of a proof of \cite[Prop.~2.82]{DaLaMu21}. 
\begin{proposition}\label{prop:coomest}
 Let $n\in {\mathbb{N}}\setminus\{0,1\}$. Let $\Omega$  be a bounded open subset of ${\mathbb{R}}^n$ of class $C^1$. Let $U$  be an open neighborhood of $\partial\Omega$  of class $C^\infty$. Let   $a\in C^\infty(\overline{U},{\mathbb{R}}^n)$ satisfy the conditions in (\ref{lem:exacu1})   for some $\vartheta\in]0,1[$ and $\tau\in]0,+\infty[$.  Let $t_1\in]0,+\infty[$ be as in Lemma \ref{lem:tubnb}. Then there exist $B\in]0,+\infty[$ and a compact subset $H$ of $\Omega$ such that
 \[
 \sup_\Omega|f|+|f:\,\Omega|_{ \omega_1	}
 \leq
 B\max\left\{\sup_H|f|, 
 \sup_H|\nabla f|, M_{t_1, \omega_1}(f)
 \right\}
 \]
 for all $f\in C^1(\Omega)$ such that 
 \begin{equation} \label{prop:coomest1}
M_{t_1, \omega_1	}(f)\equiv\sup\{
|\log|t|\,|^{-1}|\nabla f(x+ta(x))|:\,(x,t)\in (\partial\Omega)\times ]-t_1,0[ 
\}<+\infty\,,
\end{equation}
 (cf.~(\ref{omth}) for the definition of $\omega_1$). 
\end{proposition}
{\bf Proof.} Since $\Omega$ is of class $C^1$, for each point $x\in \partial \Omega$ there exists a coordinate cylinder $C(x,R_x,r_x,\delta_x)$ for $\Omega$  around $x$ and if
$r_x+\delta_x$ is less than the distance between $\partial\Omega$ and ${\mathbb{R}}^n\setminus A(t_1)$, then we have $C(x,R_x,r_x,\delta_x)\subseteq A(t_1)$.  Thus,	 we now choose $r_\ast$, $\delta_\ast\in]0,+\infty[$ such that  
\[
r_\ast+\delta_\ast<{\mathrm{dist}} (\partial\Omega,{\mathbb{R}}^n\setminus A(t_1))\,.
\]
Since we plan to invoke Lemma \ref{lem:cylcoomest}, we also assume that
\[
2\delta_\ast<\min\{\tau/2,t_1\} \,,
\]
 where $\tau$  is as  in (\ref{lem:exacu1}). 
Then, by the  lemma of  the uniform cylinders
(cf.~\cite[Lemma 2.63]{DaLaMu21}), there  exist $r\in ]0,r_\ast[$, $\delta\in]0,\delta_\ast[$, with $r<\delta$, such that if $x\in \partial\Omega$, then there exists $R_x $ in $O_n({\mathbb{R}})$ such that $C(x,R_p,r,\delta)$ is a coordinate cylinder for $\Omega$ around $x$ and the corresponding function $\gamma_x$ satisfies the condition $\nabla\gamma_x(0)=0$ and the inequalities
\begin{eqnarray*}\label{prop:coomest1a}
&&|\nabla\gamma_x(\eta)|\leq1/3\qquad\forall \eta\in \overline{{\mathbb{B}}_{n-1}(0,r)}\,,
\\ \nonumber
&& 
\sup_{x\in \partial\Omega}\|\gamma_x\|_{ C^{1}(\overline{{\mathbb{B}}_{n-1}(0,r)}) }<+\infty\,.
\end{eqnarray*}
Since $r+\delta<{\mathrm{dist}} (\partial\Omega,{\mathbb{R}}^n\setminus A(t_1))$, we have
\[
C(x,R_x,r,\delta)\cap\Omega\subseteq  A^+(t_1) = A(t_1)\cap\Omega\qquad\forall x\in \partial\Omega\,.
\]
Since $\partial\Omega$ is compact, there exists a finite family $\{x^{(j)}\}_{j=1}^m$ of points of $\partial\Omega$ such that
\[
\partial\Omega\subseteq \bigcup_{j=1}^m [x^{(j)}+R_{ x^{(j)} }^t(  {\mathbb{B}}_{n-1}(0,r/4)\times]-\delta/4,\delta/4[)]
\,,
\]
and we note that the right-hand side is an open neighborhood of $\partial\Omega$. We now set 
\begin{eqnarray*}
\lefteqn{
\mu\equiv \min_{j=1,\dots,m}\biggl\{\biggr.
  r/4,\frac{ (1-\vartheta)^{1/2}}{2\sqrt{2}({\mathrm{Lip}}(a)+1)},  \frac{t_1}{2},
}
\\ \nonumber
&&\qquad\qquad 
{\mathrm{dist}}\biggl(\biggr.
\partial\Omega, {\mathbb{R}}^n \setminus \bigcup_{j=1}^m [x^{(j)}+R_{ x^{(j)} }^t(  {\mathbb{B}}_{n-1}(0,r/4)\times]-\delta/4,\delta/4[)]
\biggl.\biggr)
\biggl.\biggr\}\,,
\end{eqnarray*}
and we choose $t_2\in]0,\mu[$. In particular, we have
\begin{eqnarray*}
\lefteqn{
A^+(t_2) \subseteq A(t_2)\subseteq\{
x\in\Omega:\,{\mathrm{dist}}(x,\partial\Omega)<\mu
\}
}
\\ \nonumber
&&\qquad
\subseteq \bigcup_{j=1}^m [x^{(j)}+R_{ x^{(j)} }^t(  {\mathbb{B}}_{n-1}(0,r/4)\times]-\delta/4,\delta/4[)]\,.
\end{eqnarray*}
Then Lemma \ref{lem:cylcoomest} implies that there exists $B_j\in]0,+\infty[$ such that
\begin{equation}\label{prop:coomest2}
|f:\,[x^{(j)}+R_{ x^{(j)} }^t(  {\mathbb{B}}_{n-1}(0,r/4)\times]-\delta/4,\delta/4[)]\cap A^+(t_2)|_{\omega_1 }\leq B_jM_{t_1, \omega_1	}(f)  
\end{equation}
for all $j\in\{1,\dots,m\}$ and for all $f\in C^1(\Omega)$ such that  $M_{t_1, \omega_1	}(f)<+\infty$. Now let $\Omega_1$ be an open subset of class $C^\infty$ of $\Omega$  such that
\[
\overline{\Omega\setminus A^+(t_2)   }\subseteq\Omega_1\subseteq \overline{\Omega_1}\subseteq \Omega
\]
(cf.~\textit{e.g.}, \cite[Lemma 2.70]{DaLaMu21}). Since $\Omega_1$ is of class $C^1$, $C^1(\overline{\Omega_1})$ is continuously embedded into $C^{0,1}(\overline{\Omega_1})$ and $C^{0,1}(\overline{\Omega_1})$ is continuously embedded into $C^{0,\omega_1 }(\overline{\Omega_1})$. Hence, there exists $c_{\omega_1}(\Omega_1)\in ]0,+\infty[$ such that
\begin{equation}\label{prop:coomest3}
|f:\,\overline{\Omega_1}|_{\omega_1 }\leq   c_{\omega_1}(\Omega_1)\sup\left\{\sup_{\overline{\Omega_1}}|f|,
\sup_{\overline{\Omega_1}}|\nabla f| \right\}
\end{equation}  
for all $f\in C^1(\overline{\Omega_1})$ (cf. (\ref{eq:0t0om0t}), \cite[Corollary 2.57 (i), (iii)]{DaLaMu21}).  Then we take  $f\in C^1(\Omega)$ such that  $M_{t_1, \omega_1	}(f)<+\infty$ and we turn to estimate $|f:\,\Omega|_{\omega_1 }$. To do so, we observe that
\begin{equation}\label{prop:coomest4}
 \overline{\Omega}\subseteq \Omega_1\cup \bigcup_{j=1}^m \left\{  [x^{(j)}+R_{ x^{(j)} }^t(  {\mathbb{B}}_{n-1}(0,r/4)\times]-\delta/4,\delta/4[)] \cap A(t_2)\right\}\,.
\end{equation}
Let $\Lambda$ be a Lebesgue number corresponding to the open cover of $\overline{\Omega}$ in the right-hand side of (\ref{prop:coomest4}). We can clearly assume that
\[
\Lambda<\delta/12\,.
\]
If $p'$, $p''\in  \Omega$ and $|p'-p''|\leq\Lambda$, then both $p'$ and $p''$ belong to at least one of the open sets in 
the right-hand side of (\ref{prop:coomest4}) and thus inequalities (\ref{prop:coomest2}) and (\ref{prop:coomest3}) imply  that
\begin{eqnarray}\label{prop:coomest4a}
\lefteqn{
|f(p')-f(p'')|
}
\\	\nonumber
&&\leq \max\biggl\{
\tilde{B} M_{t_1, \omega_1}(f),    
c_{\omega_1}(\Omega_1)\sup_{\overline{\Omega_1}}|f|,
c_{\omega_1}(\Omega_1)\sup_{\overline{\Omega_1}}|\nabla f|  
\biggr\}\omega_1(|p'-p''|)\,,
\end{eqnarray}
where
\[
\tilde{B}\equiv\max_{j\in\{1,\dots,m\}}B_j\,.
\]
In order to estimate $|f(p')-f(p'')|$ in case $|p'-p''|>\Lambda$, we need to estimate $\sup_\Omega|f|$ (cf.  Remark \ref{rem:om4}). To do so, we note that
\[
x^{(j)}-\frac{t_2}{2}a(x^{(j)})\in A(t_2)^+\subseteq \Omega\qquad\forall j\in\{1,\dots,m\}\,,
\]
and that our assumptions $t_2<\mu<r/4$, $r<\delta$ imply that
\[
x^{(j)}-\frac{t_2}{2}a(x^{(j)})\in {\mathbb{B}}_n(x^{(j)},r/8)
\subseteq   x^{(j)}+R_{ x^{(j)} }^t(  {\mathbb{B}}_{n-1}(0,r/4)\times]-\delta/4,\delta/4[) 
\]
for all $j\in\{1,\dots,m\}$. By \cite[Lemma 2.70 of the Appendix]{DaLaMu21}, there exists  an open subset  $\Omega_2$ of class $C^\infty$ of $\Omega$ such that
\[
\overline{\Omega_1}\cup\left\{
x^{(j)}-\frac{t_2}{2}a(x^{(j)}):\,j\in\{1,\dots,m\}
\right\}
\subseteq\Omega_2\subseteq\overline{\Omega_2}\subseteq\Omega\,.
\]
If $p\in \Omega\setminus\overline{\Omega_2}$, then $p\in \Omega\setminus\overline{\Omega_1}$ and
\[
p\in A(t_2)^+\subseteq \bigcup_{j=1}^m [x^{(j)}+R_{ x^{(j)} }^t(  {\mathbb{B}}_{n-1}(0,r/4)\times]
-\delta/4,\delta/4[)]
\]
and accordingly, there exists $\tilde{j}\in\{1,\dots,m\}$ such that
\[
p\in  A(t_2)^+\cap  [x^{(\tilde{j})}+R_{ x^{(\tilde{j})} }^t(  {\mathbb{B}}_{n-1}(0,r/4)\times]
-\delta/4, \delta/4[)]\,.
\]
Since both $p$ and $x^{(\tilde{j})}-\frac{t_2}{2}a(x^{(\tilde{j})})$ belong to 
\[
A(t_2)^+\cap  [x^{(\tilde{j})}+R_{ x^{(\tilde{j})} }^t(  {\mathbb{B}}_{n-1}(0,r/4)\times]
-\delta/4, \delta/4[)] 
\]
and $x^{( \tilde{j} )}-\frac{t_2}{2}a(x^{(\tilde{j})})$ belongs to $\Omega_2$,
 we have
\begin{eqnarray*}
\lefteqn{
|f(p)|\leq |f(p)-f(x^{(\tilde{j})}-\frac{t_2}{2}a(x^{(\tilde{j})}))|+|f(x^{( \tilde{j} )}-\frac{t_2}{2}a(x^{(\tilde{j})}))|
}
\\ \nonumber
&&\qquad\qquad\qquad
\leq \tilde{B}M_{t_1, \omega_1}(f)\omega_1\left(\left|p-
x^{(\tilde{j})}  +	\frac{t_2}{2}a(x^{(\tilde{j})})
\right|\right)+\sup_{\overline{\Omega_2}}|f|\,.
\end{eqnarray*}
Since both $p$ and $x^{(\tilde{j})}-\frac{t_2}{2}a(x^{(\tilde{j})})$ belong to 
\[
x^{(\tilde{j})}+R_{ x^{(\tilde{j})} }^t(  {\mathbb{B}}_{n-1}(0,r/4)\times]
-\delta/4, \delta/4[)
\]
that has a diameter less than or equal to $2(r/4)+2(\delta/4)<\delta$, we have
\[
|f(p)|\leq \tilde{B}M_{t_1, \omega_1}(f)\omega_1(\delta)+\sup_{\overline{\Omega_2}}|f|\,.
\]
If instead $p\in \overline{\Omega_2}$, we certainly have $|f(p)|\leq \sup_{\overline{\Omega_2}}|f|$.
 Hence,   
\begin{equation}\label{prop:coomest5}
\sup_\Omega|f|
\leq 
\tilde{B}M_{t_1, \omega_1	}(f)\omega_1(\delta) +\sup_{\overline{\Omega_2}}|f|\,.
\end{equation}
 Then by combining Remark \ref{rem:om4} and inequalities (\ref{prop:coomest4a}), (\ref{prop:coomest5}), we deduce that	 
\begin{eqnarray*}
\lefteqn{
|f:\Omega|_{\omega_1 }
\leq \max\biggl\{\biggr.
\max\biggl\{
\tilde{B}M_{t_1, \omega_1	}(f),  
c_{\omega_1}(\Omega_1)\sup_{\overline{\Omega_1}}| f|,
}
\\ \nonumber
&&\qquad\qquad\qquad\qquad\qquad\qquad 
c_{\omega_1}(\Omega_1)\sup_{\overline{\Omega_1}}|\nabla f|   
\biggr\}
,
\frac{2}{\omega_1(\Lambda)}
\sup_{\Omega}|f|
\biggl.\biggr\}\,.
\end{eqnarray*}
Then by taking $H=\overline{\Omega_2}$, we conclude that $B$ as in the statement  does exist.
\hfill  $\Box$ 

\vspace{\baselineskip}

\section{An extension of a theorem of C.~Miranda}
If $n\in {\mathbb{N}}\setminus\{0\}$, $m\in {\mathbb{N}}$, $h\in {\mathbb{R}}$, $\alpha\in ]0,1]$, then  we say that a function $k$ from ${\mathbb{R}}^n\setminus\{0\}$ to ${\mathbb{C}}$ is positively homogeneous of degree $h$
provided that
\[
k(tx)=t^hk(x)\qquad\forall (t,x)\in]0,+\infty[\times ({\mathbb{R}}^n\setminus\{0\})
\]
and we set 
\[ 
{\mathcal{K}}^{m,\alpha}_h \equiv\biggl\{
k\in C^{m,\alpha}_{ {\mathrm{loc}}}({\mathbb{R}}^n\setminus\{0\}):\, k\ {\text{is\ positively\ homogeneous\ of \ degree}}\ h
\biggr\}\,,
\] 
where $C^{m,\alpha}_{ {\mathrm{loc}}}({\mathbb{R}}^n\setminus\{0\})$ denotes the set of functions of 
$C^{m}({\mathbb{R}}^n\setminus\{0\})$ whose restriction to $\overline{\Omega}$ is of class $C^{m,\alpha}(\overline{\Omega})$ for all bounded open subsets $\Omega$ of ${\mathbb{R}}^n$ such that
$\overline{\Omega}\subseteq {\mathbb{R}}^n\setminus\{0\}$
and we set
\[
\|k\|_{ {\mathcal{K}}^{m,\alpha}_h}\equiv \|k\|_{C^{m,\alpha}(\partial{\mathbb{B}}_n(0,1))}\qquad\forall k\in {\mathcal{K}}^{m,\alpha}_h\,.
\]
We can easily verify that $ \left({\mathcal{K}}^{m,\alpha}_h , \|\cdot\|_{ {\mathcal{K}}^{m,\alpha}_h}\right)$ is a Banach space and we  consider the closed subspace
\[ 
{\mathcal{K}}^{m,\alpha}_{h;o} \equiv  \left\{
k\in {\mathcal{K}}^{m,\alpha}_h:\,k\ {\mathrm{is\ odd}} 
\right\}\,,
 \] 
of ${\mathcal{K}}^{m,\alpha}_h$.   
In the following Theorem \ref{thm:mirandao01}, we present an extension of a classical result of Miranda \cite{Mi65} (see also \cite[Thm.~4.17]{DaLaMu21}). In his work, Miranda considered domains of class $C^{1,\alpha}$ and densities $\mu \in C^{0,\alpha}(\partial\Omega)$ for $\alpha \in ]0,1[$; here we address the limiting case $\alpha = 1$.
\begin{theorem}\label{thm:mirandao01}
Let $\Omega$ be a bounded open subset of ${\mathbb{R}}^n$ of class $C^{1,1}$. Then the following statements hold.
\begin{enumerate}
\item[(i)] For each $(k,\mu)\in {\mathcal{K}}^{1,1}_{-(n-1);o}\times C^{0,1}(\partial\Omega)$, the map
\[
K[k,\mu]_{|\Omega}(x)= \int_{\partial\Omega} k(x-y) \mu(y)\,d\sigma_{y}\qquad \forall x\in\Omega 
\]
can be extended to a unique  $\omega_1 $-H\"{o}lder continuous function $K[k,\mu]^+$ on $\overline{\Omega}$. Moreover, the map 
from ${\mathcal{K}}^{1,1}_{-(n-1);o}\times C^{0,1}(\partial\Omega)$ to $C^{0,\omega_1 }(\overline{\Omega})$  taking	   $(k,\mu)$ to $K[k,\mu]^+$   is bilinear and continuous.
\item[(ii)]  For each   $(k,\mu)\in {\mathcal{K}}^{1,1}_{-(n-1);o}\times C^{0,1}(\partial\Omega)$  the map
\[
K[k,\mu]_{| {\mathbb{R}}^n\setminus\overline{\Omega} }(x)= \int_{\partial\Omega} k(x-y) \mu(y)\,d\sigma_{y}\qquad \forall x\in
 {\mathbb{R}}^n\setminus\overline{\Omega}\,,
\]
can be extended to a unique continuous function $K[k,\mu]^-$ on ${\mathbb{R}}^n\setminus \Omega$. Moreover, 
if $r\in]0,+\infty[$  and $\overline{\Omega}\subseteq {\mathbb{B}}_n(0,r)$, then the restriction $K[k,\mu]^-_{|\overline{{\mathbb{B}}_n(0,r)}\setminus\Omega}$ is $\omega_1 $-H\"{o}lder continuous and  the map 
from ${\mathcal{K}}^{1,1}_{-(n-1);o}\times C^{0,1}(\partial\Omega)$ to the space    $C^{0,\omega_1 }( \overline{{\mathbb{B}}_n(0,r)}\setminus\Omega)$
 taking 	$(k,\mu)$ to $K[k,\mu]^-_{|\overline{{\mathbb{B}}_n(0,r)}\setminus\Omega}$   is bilinear and continuous.
\end{enumerate}
(Cf.~(\ref{omth}) for the definition of $\omega_1$).
\end{theorem}
{\bf Proof.} We first prove statement (i). The  classical theorem of continuity for integrals depending on a parameter implies that the integral $\int_{\partial\Omega} k(x-y) \mu(y)\,d\sigma_{y}$ is continuous in $x\in{\mathbb{R}}^n\setminus\partial\Omega$.  Since $\omega_1 $-H\"{o}lder continuous functions are uniformly continuous and uniformly continuous functions in $\Omega$ can be uniquely extended to the closure $\overline{\Omega}$, it suffices to show that there exists $C>0$ such that
\[
\sup_{x\in\Omega}|K[k,\mu](x)|+ |K[k,\mu]:\,\Omega|_{\omega_1 }\leq C\|k\|_{{\mathcal{K}}^{1,1}_{-(n-1)} }\|\mu\|_{C^{0,1}(\partial\Omega)}
\]
for all $(k,\mu)\in {\mathcal{K}}^{1,1}_{-(n-1);o} \times C^{0,1}(\partial\Omega)$. 
Let $\vartheta\in]0,1[$.  By Lemma \ref{lem:exacu} there exist an open neighborhood $U$ of class $C^\infty$ of $\partial\Omega$ and $a\in C^\infty(\overline{U},{\mathbb{R}}^n)$ such that
 conditions  (\ref{lem:exacu1}) hold true for some  $\tau\in]0,+\infty[$.  Then, by the sufficient condition for H\"older continuity of Proposition \ref{prop:coomest},  it suffices to show that
\begin{enumerate}
\item[(j)] If $H$ is a compact subset of $\Omega$, then there exists $C_H>0$ such that
\[
\sup_{x\in H}|K[k,\mu](x)|+ \sup_{x\in H}|\nabla K[k,\mu](x)|\leq C_H\|k\|_{{\mathcal{K}}^{1,1}_{-(n-1)}  }\|\mu\|_{C^{0,1}(\partial\Omega)}
\]
for all $(k,\mu)\in{\mathcal{K}}^{1,1}_{-(n-1);o} \times C^{0,1}(\partial\Omega)$.
\item[(jj)]  If $t_1\in]0,+\infty[$ is as in Lemma \ref{lem:tubnb}, then there exist $t'\in]0,t_1[$ and $C_2>0$ such that
\begin{equation}\label{thm:mirandao1}
M_{t', \omega_1}(K[k,\mu]_{|\Omega})\leq C_2\|k\|_{{\mathcal{K}}^{1,1}_{-(n-1)} }\|\mu\|_{C^{0,1}(\partial\Omega)}
\end{equation}
for all $(k,\mu)\in {\mathcal{K}}^{1,1}_{-(n-1);o} \times C^{0,1}(\partial\Omega)$ (see (\ref{prop:coomest1}) for the definition of $M_{t', \omega_1}(\cdot)$). 
  \end{enumerate}
  Indeed, $t'$ would also satisfy the conditions (i), (ii) of Lemma \ref{lem:tubnb}.	 Now, the validity of (j) follows by the classical theorem of continuity for integrals depending on a parameter. Indeed, the distance from $H$ to $\partial\Omega$ is positive.   Thus it remains to prove statement (jj). In order to estimate $\nabla K[k,\mu](x+ta(x))$ for $x\in\partial\Omega$ and $t\in]-t_1,0[$  we fix $j\in\{1,\dots,n\}$ and we estimate
\[
\partial_{x_j}K[k,\mu](x+ta(x))=\int_{\partial\Omega}\partial_{x_j}k(x+ta(x)-y)\mu(y)\,d\sigma_y
\]
when $(t,x)\in ]-t_1,0[\times\partial\Omega$ for a perhaps smaller $t_1\in]0,1/e[$ and when $(k,\mu)\in {\mathcal{K}}^{1,1}_{-(n-1);o}\times C^{0,1}(\partial\Omega)$. By Lemma \ref{lem:bebelbd}   and by the properties of $a$ in  (\ref{lem:exacu1}), we have
\begin{equation}\label{thm:mirandao4}
|x-y+ta(x)|\geq (1-\vartheta)^{1/2} (|x-y|^2+|t|^2)^{1/2}
\end{equation}
for all $x$, $y\in \partial\Omega$ such that $|x-y|<\tau$ and for all $t\in ]-t_1,0[$. Next we note that if  $\mu\in C^{0,\alpha}(\partial\Omega)$, then 
\begin{eqnarray}\label{thm:mirandao5}
\lefteqn{
\int_{\partial\Omega}\partial_{x_j}k(x+ta(x)-y)\mu(y)\,d\sigma_y
}
\\ \nonumber
&&\qquad
=
\int_{\partial\Omega}\partial_{x_j}k(x+ta(x)-y)(\mu(y)-\mu(x))\,d\sigma_y
\\ \nonumber
&&\qquad\quad
+\mu(x)\int_{\partial\Omega}\partial_{x_j}k(x+ta(x)-y)\,d\sigma_y\qquad\forall x\in\partial\Omega\,
\end{eqnarray}
for all $t\in]-t_1,0[$. In order to make the proof more readable,  we now find convenient to split it 
into two parts, in which we consider separately the first and the second integral in the right-hand side of (\ref{thm:mirandao5}) and
into a series of intermediate steps. 

\noindent
{\bf Part 1.} Here we estimate  the first integral in the right-hand side of (\ref{thm:mirandao5}). To do so, we split the domain of integration $\partial\Omega$ into the two parts
\[
(\partial\Omega)\setminus {\mathbb{B}}_n(x,\tau)\,,
\qquad
(\partial\Omega) \cap  {\mathbb{B}}_n(x,\tau)\,.
\]
Since $\partial_{x_j}k$ is positively homogeneous of degree $-n$, by arguing  so as for equation (4.30) of	 Part 1 of the proof of \cite[Theorem 4.17, p.~141  with $\alpha=1$]{DaLaMu21}, we have
\begin{eqnarray}\label{thm:mirandao6}
\lefteqn{
\left|\int_{(\partial\Omega)\setminus {\mathbb{B}}_n(x,\tau)}\partial_{x_j}k(x+ta(x)-y)(\mu(y)-\mu(x))\,d\sigma_y
\right|
}
\\ \nonumber
&&\qquad
\leq \| \partial_{x_j}k\|_{C^0(\partial{\mathbb{B}}_n(0,1))}|\mu:\, \partial\Omega|_1
m_{n-1}(\partial\Omega){\mathrm{diam}}(\partial\Omega) (\tau/2)^{-n}
\end{eqnarray}
for all $t\in ] -\min\{ t_1,\tau/2\}, 0  [$. 

\noindent$\bullet$ We now turn   to the integral on $(\partial\Omega) \cap  {\mathbb{B}}_n(x,\tau)$. 
As  in equation (4.31) of \cite[Theorem 4.17, p.~141  	with $\alpha=1$]{DaLaMu21}, by inequality (\ref{thm:mirandao4}), we have
\begin{eqnarray}\label{thm:mirandao7}
\lefteqn{
\left|\int_{(\partial\Omega)\cap {\mathbb{B}}_n(x,\tau)} \partial_{x_j}k(x+ta(x)-y)(\mu(y)-\mu(x))\,d\sigma_y
\right|
}
\\ \nonumber
&&\qquad
\leq \frac{\| \partial_{x_j}k\|_{C^0(\partial{\mathbb{B}}_n(0,1))}|\mu:\, \partial\Omega|_1}{(1-\vartheta)^{n/2}}
\int_{ (\partial\Omega)\cap {\mathbb{B}}_n(x,\tau) } \frac{|x-y|\,d\sigma_y}{ (|x-y|^2+t^2)^{n/2}}\,.
\end{eqnarray}
\noindent$\bullet$ We now estimate the last integral in the 
right-hand side of (\ref{thm:mirandao7}) by considering  separately the case where
\[
|x-y|\geq |t|  
\]
and the case where
\[
|x-y| < |t|  
\]
for $t\in ] -\min\{ t_1,\tau/2\}, 0  [$. Then we have 
\begin{eqnarray}\label{thm:mirandao8}
\lefteqn{
\int_{ (\partial\Omega)\cap {\mathbb{B}}_n(x,\tau) } 
\frac{|x-y|\,d\sigma_y}{ (|x-y|^2+t^2)^{n/2}}
}
\\ \nonumber
&&\qquad
\leq
\int_{ \{y\in (\partial\Omega)\cap {\mathbb{B}}_n(x,\tau):\, |x-y|\geq |t|\}
} \frac{|x-y|\,d\sigma_y}{ (|x-y|^2+t^2)^{n/2}}
\\ \nonumber
&&\qquad\quad
+
\int_{\{y\in (\partial\Omega)\cap {\mathbb{B}}_n(x,\tau):\, |x-y|< |t|\} } \frac{|x-y|\,d\sigma_y}{ (|x-y|^2+t^2)^{n/2}}\,.
\end{eqnarray}
We note that Lemma \ref{lem:comin} (iv) implies that
\begin{eqnarray}\label{thm:mirandao9}
\lefteqn{
\int_{ \{y\in (\partial\Omega)\cap {\mathbb{B}}_n(x,\tau):\, |x-y|\geq |t|\}
} \frac{|x-y|\,d\sigma_y}{ (|x-y|^2+t^2)^{n/2}}
}
\\ \nonumber
&&\qquad
\leq 
\int_{ \{y\in (\partial\Omega)\cap {\mathbb{B}}_n(x,\tau):\, |x-y|\geq |t|\}
} |x-y|^{1-n}\,d\sigma_y
\\ \nonumber
&&\qquad
\leq
c_{\Omega}^{(iv)}		|\log|t|\,| 
\end{eqnarray}
for all $t\in ] -\min\{ t_1,\tau/2\}, 0  [$ and that Lemma \ref{lem:comin} (ii) implies that 
\begin{eqnarray}\label{thm:mirandao10}
\lefteqn{
\int_{\{y\in (\partial\Omega)\cap {\mathbb{B}}_n(x,\tau):\, |x-y|< |t|\} } 
\frac{|x-y|\,d\sigma_y}{ (|x-y|^2+t^2)^{n/2}}
}
\\ \nonumber
&&\qquad
\leq |t|^{-n}\int_{\{y\in (\partial\Omega)\cap {\mathbb{B}}_n(x,\tau):\, |x-y|< |t|\} } 
\frac{d\sigma_y}{|x-y|^{-1}}
\\ \nonumber
&&\qquad
\leq |t|^{-n}
c_{\Omega,-1}''|t|^{(n-1)-(-1)} 
=c_{\Omega,-1}''  
\end{eqnarray}
for all $t\in ] -\min\{ t_1,\tau/2\}, 0  [$.

\noindent$\bullet$ Finally,  combining inequalities (\ref{thm:mirandao6})--(\ref{thm:mirandao10}), we obtain 
\begin{eqnarray}\label{thm:mirandao11}
\lefteqn{
\left|\int_{\partial\Omega} \partial_{x_j}k(x+ta(x)-y)(\mu(y)-\mu(x))\,d\sigma_y
\right|
}
\\ \nonumber
&&\qquad
\leq \| \partial_{x_j}k\|_{C^0(\partial{\mathbb{B}}_n(0,1))}|\mu:\, \partial\Omega|_1
\biggl\{\biggr.
m_{n-1}(\partial\Omega){\mathrm{diam}}(\partial\Omega) (\tau/2)^{-n}
\\ \nonumber
&&+
\frac{
c_{\Omega }^{(iv)}  +c_{\Omega,-1}'' 
}{ ( 1-\vartheta)^{n/2}} 
\biggl.\biggr\}  |\log|t|\,| 		
\end{eqnarray}
for all $t\in ] -\min\{ t_1,\tau/2\}, 0  [$ and $x\in\partial\Omega$. 

\noindent {\bf Part 2.} We now turn    to 	  the second integral in the right-hand side of (\ref{thm:mirandao5}). By following the lines of the Part 2 of the proof of \cite[Theorem 4.17, p.~143]{DaLaMu21} and replacing $\alpha$ by $1$ with the sole exception of  the term
\[
\frac{|t|^{\alpha-1}}{1-\alpha}
\]
in equation (4.51) of the proof  of \cite[Theorem 4.17, p.~148]{DaLaMu21}  that has to be replaced by
\[
|\log |t||\, ,
\]
one can show the existence of $C''>0$ such that
\begin{equation}\label{thm:mirandao12}
\left|
\int_{\partial\Omega}\partial_{x_j}k(x+ta(x)-y)\,d\sigma_y
\right|
\leq C''  \| \partial_{x_j}k\|_{C^{0,1}(\partial{\mathbb{B}}_n(0,1))}  |\log|t||
\end{equation}
for all $
 t\in ] -\min\{ t_1,\tau/2,r/2 \}, 0 [=] -\min\{ t_1, r/2 \}, 0 [
$, $x\in\partial\Omega$  and $k\in {\mathcal{K}}^{1,1}_{-(n-1);o}$. Then by combining equality  (\ref{thm:mirandao5}) with inequalities (\ref{thm:mirandao11}) and (\ref{thm:mirandao12}), we conclude 
 that there exists a constant $C_2>0$  as in (\ref{thm:mirandao1}) 
 for all  $t\in ]-\min\{ t_1, r/2 \}, 0 [$  and that accordingly statement (jj) holds true (see also \cite[Lemma 4.13]{DaLaMu21}). Thus the proof of (i)  is complete.

 Then statement (ii) can be deduced applying statement (i) to the open bounded set $\mathbb{B}_n(0,r)\setminus\overline{\Omega}$. Indeed the map from $C^{0,1} (\partial\Omega)$ to $C^{0,1}(\partial(\mathbb{B}_n(0,r)\setminus\overline{\Omega}))$ that  takes $\mu$ to
\[
\tilde\mu(x)\equiv
\begin{cases}
\mu(x)&\text{if }x\in\partial\Omega\,,\\
0&\text{if }x\in\partial\mathbb{B}_n(0,r)
\end{cases}
\]
is linear and continuous and we have
\[ K[k,\mu]^-_{|\overline{\mathbb{B}_n(0,r)}\setminus\Omega}=K[k,\tilde\mu]^+ \]
for all $(k,\mu)\in  {\mathcal{K}}^{1,1}_{-(n-1);o}
\times C^{0,1}(\partial\Omega)$.

\hfill  $\Box$ 

\vspace{\baselineskip}

  \noindent
{\bf Statements and Declarations}\\

 \noindent
{\bf Competing interests:} This paper does not have any  conflict of interest or competing interest.

 \noindent
{\bf Acknowledgement.} 

%The author  is indebted  to Prof.~Paolo Musolino for a number of comments on the paper.   

The authors are members of the ``Gruppo Nazionale per l'Analisi Matematica, la Probabilit\`a e le loro Applicazioni'' (GNAMPA) of the ``Istituto Nazionale di Alta Matematica'' (INdAM). The authors  acknowledge the support  of the Project funded by the European Union – Next Generation EU under the National Recovery and Resilience Plan (NRRP), Mission 4 Component 2 Investment 1.1 - Call for tender PRIN 2022 No. 104 of February, 2 2022 of Italian Ministry of University and Research; Project code: 2022SENJZ3 (subject area: PE - Physical Sciences and Engineering) ``Perturbation problems and asymptotics for elliptic differential equations: variational and potential theoretic methods''.	 M. Lanza de Cristoforis and P. Musolino also acknowledge the support from EU through the H2020-MSCA-RISE-2020 project EffectFact, Grant agreement ID: 101008140.

 \end{document}